\documentclass[14pt]{amsart}

\usepackage{amsmath,amssymb,amsfonts,enumerate,amsthm, amscd}
\usepackage[english]{babel}
\usepackage[colorlinks=true]{hyperref}
\usepackage{paralist}
\usepackage{verbatim}
\usepackage{algorithm}
\usepackage{algorithmicx}
\usepackage{algpseudocode}
\usepackage{color}
\usepackage{algpseudocode}
\usepackage{xcolor} 
\usepackage{tikz}

\theoremstyle{plain}
\newtheorem{thm}{\bf Theorem}[section]

\newtheorem{prop}[thm]{\bf Proposition}
\newtheorem{defn-prop}[thm]{\bf Definition and Proposition}
\newtheorem{lem}[thm]{\bf Lemma}
\newtheorem{cor}[thm]{\bf Corollary}

\theoremstyle{definition}
\newtheorem{defn}[thm]{\bf Definition}
\theoremstyle{remark}
\newtheorem{rem}[thm]{\bf Remark}
\newtheorem{exam}[thm]{\bf Example}

\theoremstyle{example}

\def \h{\mathrm{h}}

\def \P{\mathbb{P}}

\def \F{\mathbb{F}}
\def \Z{\mathbb{Z}}
\def \RR{\mathbb{R}}
\def \Q{\mathbb{Q}}
\def \N{\mathbb{N}}

\def \L{\mathcal{L}}
\def \O{\mathcal{O}}

\def \C{\mathcal{C}}

\def \FF{\mathcal{F}}

\def \Pic{\mathrm{Pic}}

\def \rank{\mathrm{rank}}

\def \Supp{\mathrm{Supp}}

\def \ev{\mathrm{ev}}
\def \im{\mathrm{Im}}

\def \Div{\mathrm{Div}}

\begin{document}

\title[Hierarchical filtrations of line bundles and optimal AG codes]{Hierarchical filtrations of line bundles and optimal algebraic geometry codes}
\author{Rahim Rahmati-asghar}

\keywords{Line bundles, Hierarchical filtration, AG codes}

\makeatletter
%\@namedef{subjclassname@2020}{%
  %\textup{2020} Mathematics Subject Classification}
%\makeatother

\subjclass[2020]{primary; 14G50, 94B27 Secondary; 14C20, 14Q05}
%\thanks{The author would like to sincerely thank Professor Pierre Deligne for his insightful comments and critical feedback on the definition of hierarchical depth, which is the fundamental concept of this paper.}

\date{\today}

\begin{abstract}
We introduce \emph{hierarchical depth}, a new invariant of line bundles and divisors, defined via maximal chains of effective sub-line bundles. This notion gives rise to \emph{hierarchical filtrations}, refining the structure of the Picard group and providing new insights into the geometry of algebraic surfaces. We establish fundamental properties of hierarchical depth, derive inequalities through intersection theory and the Hodge index theorem, and characterize filtrations that are Hodge-tight.

Using this framework, we develop a theory of \emph{hierarchical algebraic geometry codes}, constructed from evaluation spaces along these filtrations. This approach produces nested families of codes with controlled growth of parameters and identifies an optimal intermediate code maximizing a utility function balancing rate and minimum distance. Hierarchical depth thus provides a systematic method to construct AG codes with favorable asymptotic behavior, linking geometric and coding-theoretic perspectives.

Our results establish new connections between line bundle theory, surface geometry, and coding theory, and suggest applications to generalized Goppa codes and higher-dimensional evaluation codes.
\end{abstract}

\maketitle

\section*{Introduction}

The theory of line bundles and divisors on algebraic varieties lies at the  intersection of algebraic geometry, number theory, and information theory.  Filtrations of vector bundles play a central role in understanding stability conditions, cohomological behavior, and geometric invariants (\cite{HaNa,AtBo}). In coding theory, since the foundational work of Goppa \cite{Go}, line bundles on curves have provided the framework for constructing \emph{algebraic geometry (AG) codes}, leading to significant progress on the construction of long linear codes with good asymptotic properties \cite{St,TsVlZi}. 

In this paper we introduce a new invariant of line bundles and divisors,  which we call the \emph{hierarchical depth}. This invariant arises from maximal chains of effective sub-line bundles and leads naturally to the notion of \emph{hierarchical filtrations}. While classical invariants such as degree or dimension describe line bundles globally, hierarchical depth captures the internal structure of a line bundle in terms of successive reductions. Our approach provides a systematic way of encoding information about positivity and effectivity in divisor theory. 

We establish several foundational properties of hierarchical depth. In particular, we show that every effective line bundle admits a hierarchical filtration of bounded length, and we derive inequalities controlling the depth using intersection theory and the Hodge index theorem \cite{Ha,Be}. For surfaces, we introduce the notion of \emph{Hodge-tight filtrations}, in which each step in the filtration saturates the Hodge index bound. We illustrate these constructions with examples from projective spaces, elliptic curves, and Hermitian curves, which demonstrate both the flexibility and the constraints of the theory. 

A central motivation for our work comes from applications to coding theory. By evaluating global sections of line bundles along a hierarchical filtration, we obtain nested sequences of AG codes, which we call \emph{hierarchical AG codes}. This viewpoint generalizes the classical construction of Goppa codes and provides new tools for analyzing code parameters. In particular, we show that within such a hierarchical family there exists a \emph{distinguished middle layer} code that optimally balances rate and minimum distance, as measured by a natural utility function. This construction sheds light on the geometry underlying the trade-offs between code parameters and suggests new approaches to asymptotic bounds \cite{TsVlZi,St}. 

While nested families of AG codes have appeared previously in the literatureunder the form of flag codes \cite{BrAmDuHo}, the present work differs in that the filtrations are canonically induced by hierarchical depth and are analyzed through explicit optimization criteria on curves and surfaces.

The novelty of this work lies not in the existence of nested AG code chains per se, which are classical, but in introducing hierarchical depth as a new geometric invariant and using it to define canonical filtrations whose optimization yields explicit best-code selection principles. In particular, hierarchical depth is not merely nested inclusion. It is a new invariant measuring minimal filtration length needed to generate a divisor/line bundle.

Our results point to several avenues for further research. From the geometric side, hierarchical depth offers a new invariant for the classification of line bundles, raising questions about its behavior on higher-dimensional varieties and its relation to stability conditions. From the coding-theoretic side, hierarchical codes suggest a framework for refining classical bounds, and designing evaluation codes over higher-dimensional varieties.

\medskip
\noindent
\textbf{Outline.} In Section \ref{Sec:Hierar} we introduce hierarchical depth  and filtrations and establish their basic properties on curves and surfaces. In the sequel, we exploit the Hodge index theorem to obtain tight inequalities and structural results on hierarchical depth. In Section \ref{AG codes} we develop the theory of hierarchical AG codes and prove the existence of optimal middle-layer codes, and compare these constructions with classical AG codes. Throughout the paper we try to present more examples to clarify the constructions.

\section*{Acknowledgments}

The author would like to express his deepest gratitude to Professor Angelo López for his generous support, valuable guidance, and continuous encouragement throughout the development of this work. His insights and suggestions were instrumental in shaping the ideas presented in this paper.

The author also sincerely thanks Professor Pierre Deligne for his insightful comments and critical feedback on the definition of hierarchical depth, which forms the central concept of this work. In addition, the author is grateful to Professor Carles Padró for kindly reviewing an early draft of the paper and offering valuable suggestions that contributed to its improvement.

\section{Hierarchical depth and filtrations}
\label{Sec:Hierar}

\begin{defn}\label{Def:hierar}
Let $X$ be a smooth projective variety over a field $\F$ and let $L$ be a line bundle on $X$.

A \emph{hierarchical filtration} of $L$ is a finite chain of inclusions of coherent subsheaves
$$\FF_L:\O_X = \L_0 \subset \L_1 \subset \cdots \subset \L_h = L$$
such that for each $i=1,\dots,h$, there exists a nonzero effective Cartier divisor $E_i$ on $X$ satisfying:
$$\L_i\cong\L_{i-1}\otimes\O_X(E_i),$$
or equivalently, the quotient $\L_i / \L_{i-1}$ is a nonzero torsion sheaf supported on an effective Cartier divisor, i.e., there exists a nonzero section
$$s_i\in H^0(S,\L_i\otimes\L^{-1}_{i-1})$$
whose vanishing locus defines a nonzero effective divisor.

The maximal length of all finite hierarchical filtrations of $L$ is defined as \textit{hierarchical depth} of $L$ and denoted by $h(L)$. We will imply to the existence of a filtration of maximal length in Proposition \ref{prop:max-filtration}.

We define $h(\O_X)=0$ and if no hierarchical filtration exists for $L$, then we set $h(L)=-\infty$.
\end{defn}

Note that $\L_{i-1}\subset\L_i$ is interpreted as inclusion of coherent sheaves or line subbundles. Meanwhile, the filtration provides a discrete, layered way of analyzing how sections of line bundles build up via effective divisors. It encodes geometric information into an algebraic length invariant.

\begin{prop}\label{prop:max-filtration}
Let $X$ be a smooth projective variety over a field, and let $L$ be a line bundle on $X$.
If there exists at least one hierarchical filtration of $L$, then the set of possible
lengths of such filtrations is bounded above and hence admits a maximum. 
In particular, the hierarchical depth $h(L)$ is finite, and there exists a hierarchical
filtration of $L$ of maximal length.
\end{prop}

\begin{proof}
Fix an ample line bundle $\O_X(1)$ on $X$, and let $H=c_1(\O_X(1))\in \Pic(X)$. 
Let 
$$
\O_X=\L_0 \subset \L_1 \subset \cdots \subset \L_h=L
$$
be a hierarchical filtration of $L$. 
By definition, for each $i$ there exists a nonzero effective Cartier divisor $E_i$ on $X$ 
such that 
$$
\L_i \cong \L_{i-1}\otimes \O_X(E_i).
$$
Equivalently, there is a nonzero section 
$s_i \in H^0\!\big(X,\L_i\otimes\L_{i-1}^{-1}\big)$ 
whose divisor of zeros is precisely $E_i$. 

Multiplying the sections $s_1,\dots,s_h$ yields a nonzero section 
$$
s := s_1 \cdots s_h \in H^0(X,L),
$$
and its zero divisor satisfies
$$
\operatorname{div}(s) = E_1 + \cdots + E_h.
$$
Now intersecting with $H^{\dim X - 1}$ gives
$$
\sum_{i=1}^h (E_i \cdot H^{\dim X - 1})
= \operatorname{div}(s)\cdot H^{\dim X - 1}.
$$
Since $\operatorname{div}(s)$ is linearly equivalent to $c_1(L)$, the right-hand side depends
only on $L$, not on the chosen filtration. Set
$$
N := c_1(L)\cdot H^{\dim X - 1} \in \mathbb{Z}_{\ge 0}.
$$
Because $X$ is projective and $H$ is ample, one has 
$E\cdot H^{\dim X-1} > 0$ for every nonzero effective Cartier divisor $E$ on $X$. 
Hence each summand $E_i\cdot H^{\dim X-1}$ is at least $1$, so
$$
h \le \sum_{i=1}^h (E_i \cdot H^{\dim X - 1}) = N.
$$
Thus the length $h$ of any hierarchical filtration of $L$ is bounded above by $N$.
Consequently the set of all possible lengths is a nonempty finite subset of 
$\{0,1,2,\dots,N\}$, and therefore admits a maximum. 
By definition this maximum is $h(L)$, and any filtration realizing it is a 
filtration of maximal length.
\end{proof}

\begin{rem}
The integer 
$$
N = c_1(L)\cdot H^{\dim X - 1}
$$
is the degree of $L$ with respect to the polarization $H=\O_X(1)$. 
Thus Proposition \ref{prop:max-filtration} shows that the hierarchical depth of $L$
is always finite and bounded above by $\deg_H(L)$. 
\end{rem}

\begin{defn}
For a divisor $D$ on a smooth projective variety $X$, the hierarchical depth $h(D)$ is defined $h(D):=h(\O_X(D))$.
\end{defn}

Given the hierarchical filtration of line bundles
$$\L_0\subset\L_1\subset\ldots\subset L_h=\O(D),$$
each step satisfies
$$\L_i=\L_{i-1}\otimes\O(E_i)$$
where $E_i$ is a divisor such that $\O(E_i)$ has a nonzero section, meaning that $E_i$ is effective. Taking $D_i=D_{i-1}+E_i$, we obtain the sequence
$$0=D_0<D_1<\dots<D_h=D,$$
which each step satisfies $D_i-D_{i-1}=E_i$ which is effective. The inequality $D_{i-1}\leq D_i$ and $D_{i-1}< D_i$ are to be interpreted according to the notion defined in \cite{Ha}.

Note that for a line bundle $L=\O_X(D)$ on $X$, $h(L)>0$ only holds if $D$ is effective and $D>0$. Actually, because for all $i$, $D_i-D_{i-1}$ is effective, it follows that $D=D-D_0=\sum(D_i-D_{i-1})$ is effective.

\subsection{Hierarchical filtrations on curves}

Let $C$ be a smooth projective algebraic curve over a field. Because effective divisors on $C$ are of the form $D=\sum n_iP_i$ where $P_i$'s are points on $C$ and $n_i>0$, so it follows from the definition of hierarchical depth that $h(L)=\deg(D)$ where $L=\O_C(D)$. We formally address this fact in the following as a corollary of Proposition \ref{prop:max-filtration}.

\begin{cor}\label{cor:curve-degree}
Let $C$ be a smooth projective curve over a field and let $L$ be a line bundle on $C$. If $L$ admits a hierarchical filtration then every hierarchical filtration of $L$ has length at most $N=\deg(L)$. Moreover, if
$H^0(C,L)\neq 0$ then $h(L)=\deg(L)$ and else there is no hierarchical filtration.
\end{cor}

\begin{proof}
Since $C$ is a curve, it follows from Proposition \ref{prop:max-filtration} that 
$$h(L)\leq \deg(L).$$

If $H^0(C,L)\neq0$ choose a nonzero section $s\in H^0(C,L)$; its divisor of zeros
$\operatorname{div}(s)$ is an effective divisor $D$ of degree $\deg(L)$ and one has
$L\cong\O_C(D)$. Write $D=\sum_{j=1}^{\deg(L)} P_j$ where the points $P_j$ appear
with multiplicity (choose an ordering of the points repeating according to multiplicity).
For $i=0,\dots,\deg(L)$ set $D_i=\sum_{j=1}^i P_j$ (so $D_0=0$ and $D_{\deg(L)}=D$)
and define
$$
\L_i := \O_C(D_i).
$$
Each inclusion $\L_{i-1}\subset\L_i$ is given by tensoring with $\O_C(P_i)$, equivalently
the quotient $\L_i/\L_{i-1}$ is a nonzero torsion sheaf supported at the point $P_i$.
Thus
$$
\O_C=\L_0\subset\L_1\subset\cdots\subset\L_{\deg(L)}=L
$$
is a hierarchical filtration of $L$ of length $\deg(L)$. Combined with the inequality
$h(L)\le\deg(L)$ this shows $h(L)=\deg(L)$ when $H^0(C,L)\neq0$.

If $H^0(C,L)=0$ then by our convention $h(L)=-\infty$ and there is no hierarchical filtration.
\end{proof}

The following corollary is a straightforward result of Riemann-Roch theorem and Corollary \ref{cor:curve-degree}.

For a variety $X$ and a line bundle on it, we set $\h^i(L):=\dim H^i(X,L)$.

\begin{cor}
Let $C$ be a smooth projective curve over a field and let $L$ be a line bundle on $C$. When $C$ is of genus $g$ and $L=\O_C(D)$ a line bundle on $C$, then
\begin{equation}
h(L)\leq\h^0(C,L)+g-1.
\label{Eq.h(L)leq}
\end{equation}
In particular, the equality holds when $\deg(D)>2g-2$.
\end{cor}

\begin{proof}
Note that $\deg(D)>2g-2$ gives $\h^1(C,L)=0$ and it follows from Riemann-Roch that $\ell (D)=\deg(D)-g+1$.
\end{proof}

\subsection{Hierarchical filtrations on surfaces}

Let $X=S$ be a smooth projective surface and let $L$ be a line bundle on $S$.
Fix an ample polarization $\O_S(1)$ and write $H:=c_1(\O_S(1))\in\Pic(S)$.
Recall from Proposition \ref{prop:max-filtration} that
$$
N  =  c_1(L)\cdot H^{\dim S-1}  =  c_1(L)\cdot H \in\mathbb{Z}_{\ge 0},
$$
and that for any hierarchical filtration of $L$ of length $h$ one has $h\le N$.

The intersection number $c_1(L)\cdot H$ is the degree of $L$ with respect to the
polarization $\O_S(1)$; equivalently it is the intersection product of the divisor
class of $L$ with the ample class $H$. Concretely, if $s\in H^0(S,L)$ is a
nonzero section with divisor of zeros $D=\operatorname{div}(s)$ (an effective Cartier
divisor), then
$$
N  =  D\cdot H  =  \sum_{i} (E_i\cdot H)
$$
for any decomposition $D=\sum_i E_i$ into effective Cartier divisors. Since each
term $E_i\cdot H$ is a positive integer, the inequality $h\le N$ follows as in
Proposition \ref{prop:max-filtration}.

\begin{rem}[When the bound is attained]
The inequality $h(L)\le N$ is in general strict or an equality depending on the
existence of a section $s\in H^0(S,L)$ whose zero divisor decomposes into many
effective Cartier summands of small $H$-degree. More precisely:
\begin{itemize}
  \item If there exists a section $s$ with $\operatorname{div}(s)=\sum_{j=1}^m E_j$
    where each $E_j$ satisfies $E_j\cdot H =1$, then $m = D\cdot H = N$ and the
    filtration
    $$
      \O_S = \O_S \bigl(0\bigr)
      \subset \O_S(E_1)
      \subset \O_S(E_1+E_2)
      \subset \cdots
      \subset \O_S\!\bigl(\sum_{j=1}^m E_j\bigr) = L
    $$
    is a hierarchical filtration of length $m=N$. Thus in this situation the bound
    is attained and $h(L)=N$.
  \item Conversely, if every nonzero section $s\in H^0(S,L)$ has divisor $D=\operatorname{div}(s)$
    which cannot be written as a sum of $N$ (or more) effective Cartier divisors each
    of $H$-degree $1$, then $h(L)<N$. In particular, for a general section whose
    divisor is irreducible one necessarily has $h(L)=1$ even when $N$ is large.
\end{itemize}
\end{rem}

In the sequel we investigate $h(L)$ from cohomological dimension point of view.

\begin{thm}
\label{thm:hierar_formula}
Let $S$ be a smooth projective surface over an algebraically closed field $\F$, and let 
$$0=D_0< D_1<\ldots< D_h=D\quad \text{with} \quad E_i:=D_i-D_{i-1}>0$$
be a nested sequence of effective divisors. Let $L = \O_S(D)$. Then there is an exact sequence for each step:
$$0\longrightarrow \O_S(D_{i-1})\longrightarrow\O_S(D_i)\longrightarrow\O_{E_i}(D_i) 
\longrightarrow 0.$$
This yields:
\begin{equation}
\label{Formula:hierar}
\h^0(S,L)=1+\sum_{i=1}^h \big(\h^0(E_i,\O_{E_i}(D_i))-t_i\big),
\end{equation}
where 
$$t_i:=\dim\im\big(H^0(E_i,\O_{E_i}(D_i))\to H^1(S,\O_S(D_{i-1}))\big).$$
\end{thm}
\begin{proof}
For each $i=1,\dots,h$, the standard short exact sequence
$$0\longrightarrow\O_S(D_{i-1})\longrightarrow\O_S(D_i)\longrightarrow \O_{E_i}(D_i)\longrightarrow 0.$$
of sheaves yields a long exact sequence in cohomology:

\begin{flushleft}
$0\longrightarrow H^0(S,\O_S(D_{i-1}))\longrightarrow H^0(S, \O_S(D_i))\overset{\lambda_i}{\longrightarrow}$
\end{flushleft}
\begin{flushright}
\begin{equation}
H^0(E_i, \O_{E_i}(D_i))\overset{\delta_i}{\longrightarrow} H^1(S, \O_S(D_{i-1}))
\longrightarrow \cdots .
\label{seq}
\end{equation}
\end{flushright}

By dimension count we have
$$\h^0(S, D_i) 
=\h^0(S, D_{i-1})+\h^0(E_i,\O_{E_i}(D_i))-t_i,\quad 
\text{where } t_i = \dim\im(\delta_i).$$

Iterations give us:
$$
\begin{aligned}
\h^0(S, D) &=\h^0(S, D_h) 
=\h^0(S, D_0) + \sum_{i=1}^h \Big(\h^0(E_i, \O_{E_i}(D_i)) - t_i \Big).
\end{aligned}
$$
Since $D_0=0$, we have $\h^0(S, D_0)=\h^0(S,\O_S)=1$. Hence,
$$\h^0(S,L)=1+\sum_{i=1}^h\big(\h^0(E_i,\O_{E_i}(D_i))-t_i\big).$$
This completes the proof.
\end{proof}

In a view to Formula (\ref{Formula:hierar}) and the long exact sequence in above proof we have:

\begin{align}
\h^0(E_i,\O_{E_i}(D_i))-t_i &=\dim(\ker\delta_i) \nonumber\\
 &=\dim(\im\lambda_i) \nonumber\\
 &=\h^0(S,\O_S(D_i))-\dim(\ker\lambda_i) \nonumber\\
 &=\h^0(S,\O_S(D_i))-\h^0(S,\O_S(D_{i-1})) \label{Eq:h^0-t}
\end{align}

This implies that

\begin{equation}
\label{Formula2:hierar}
\h^0(S,L)=1+\sum_{i=1}^h\big(\h^0(S,\O_S(D_i))-\h^0(S,\O_S(D_{i-1}))\big).
\end{equation}

Formula (\ref{Formula:hierar}) emphasizes the contribution at each step coming from sections on the divisors $E_i$ and possible obstructions measured by $t_i$. 
It shows how the geometry of the intermediate divisors controls the growth of sections. Each $E_i$ gives potential new sections, while the maps in \eqref{seq} may kill them if they lift nontrivially to cohomology. While, formula (\ref{Formula2:hierar}) rewrites the same count purely in terms of the differences in the global sections along the filtration.
It highlights that the hierarchical depth $h$ provides an upper bound for how many times the global sections can strictly increase.
This version is sometimes more practical in concrete calculations, because one may directly compute $\h^0$ for the intermediate line bundles, for example by using Riemann-Roch and vanishing theorems.

Both versions together show that the hierarchy bridges local contributions on the supports $E_i$ and the global structure of sections on $S$.

\begin{prop}
\label{prop:hierar_dA}
Let $X$ be a smooth projective variety such that $\Pic(X)\cong\Z\cdot A$ for some ample effective divisor $A$. Then for any line bundle $L=\O_X(dA)$, the hierarchical depth satisfies $h(L)=d.$
\end{prop}

\begin{proof}
We aim to show that the maximal length of a hierarchical filtration of $L=\O_X(dA)$ is exactly $d$.

\textbf{Step 1.}
Consider the ascending chain of line bundles
$$\O_X=\O_X(0A)\subset\O_X(A)\subset\O_X(2A)\subset\cdots\subset\O_X(dA)=L.$$
Define the filtration
$$\L_i:=\O_X(iA),\quad \text{for } i=0,1,\ldots,d.$$
Then each successive quotient satisfies
$$\L_i / \L_{i-1} \cong \O_X(iA)/\O_X((i-1)A) \cong \O(D_i),$$
where $ D_i $ is the divisor of a nonzero section in $H^0(X,\O_X(A))$, since $A$ is effective and ample. Each such quotient is supported on a Cartier divisor $D_i\sim A$, and hence the filtration satisfies the conditions of Definition~\ref{Def:hierar}. Thus, this gives a hierarchical filtration of length $d$, so $h(L)\geq d.$

\textbf{Step 2.} % Maximality of the length $d$:}
Now suppose we have any hierarchical filtration of $L$:
$$\L_0 \subsetneq \L_1 \subsetneq \cdots \subsetneq \L_h = L,$$
with each $\L_i/\L_{i-1}\cong\O(E_i)$, for some effective divisors $E_i>0$.

Then since all line bundles are of the form $\O_X(k_i A)$, we may write $\L_i\cong\O_X(k_i A)$ for some integers $0=k_0<k_1<\cdots <k_h=d$. Because each inclusion is strict and corresponds to an effective divisor, we have
$$k_i-k_{i-1}\geq 1\quad\text{for each } i.$$
Summing up
$$\sum_{i=1}^h(k_i-k_{i-1})=k_h-k_0=d.$$
Since each increment is at least 1, the number of steps satisfies $h\leq d.$

From the lower and upper bounds, we conclude $h(L)=d.$
\end{proof}

\begin{rem}
The condition $\Pic(X)\cong \Z\cdot A$ is restrictive, but holds for projective spaces, smooth complete intersections of sufficiently low degree, and certain hypersurfaces.
\end{rem}

\begin{exam}
Let $S=\P^2$, the complex projective plane. The canonical bundle is $K_S=\O_{\P^2}(-3)$, and $\Pic(\P^2)\cong\Z\cdot H$, where $H$ denotes the class of a line.

Let $L=\O_{\P^2}(d)$. Fix $d=2$, so $L=\O(2)$. We examine the hierarchical depth $h(L)$, as defined via a filtration of line bundles
$$\O=\L_0\subsetneq\L_1\subsetneq\cdots\subsetneq \L_h=\O(2),$$
such that each quotient $ \L_i / \L_{i-1} \cong \O(D_i) $, with $ D_i > 0 $ effective Cartier divisors.

Applying the Riemann-Roch formula gives
$$h^0(\O(2))=\frac{1}{2}(d^2 + 3d)+1=\frac{1}{2}(4+6)+1=6.$$

In particular, for $d=2$, Proposition \ref{prop:hierar_dA} implies that $h(\O(2))=2.$

Consequently,
$$h^0(\O(2))=6\geq h(\O(2))=2.$$
\end{exam}

\subsection{A high bound for hierarchical depth of Hodge-tight filtrations}

Let $S$ be a smooth projective surface over an algebraically closed field,
and let $H$ be an ample divisor on $S$ or $H^2>0$. Then as a known result of Hodge-index theorem \cite[Theorem V.1.9]{Ha}, for every divisor $D\in\Div(S)$ one has

\begin{equation}\label{eq:hodge}
(H \cdot D)^2 \ge H^2 D^2.
\end{equation}
See \cite[Ex. V.1.9]{Ha}.

Let $NS(S)$ denote the Néron-Severi group of $S$ (divisor classes modulo algebraic equivalence \cite[p. 367]{Ha}) and write
$$NS(S)_{\RR}=NS(S)\otimes_{\Z}\RR,\qquad NS(S)_{\Q}=NS(S)\otimes_{\Z}\Q.$$
The intersection product on divisors induces a nondegenerate symmetric bilinear form
$$(\ ,\ ) : NS(S)_{\RR}\times NS(S)_{\RR}\to\RR,\qquad (X,Y)=X\cdot Y.$$
By the Hodge index theorem the intersection form on $NS(S)_{\RR}$ has signature $(1,\rho-1)$ where $\rho=\rank NS(S)$, and the $1$-dimensional positive direction may be taken to be the ray spanned by the ample class $H$. In particular the orthogonal complement
$$H^\perp=\{x\in NS(S)_{\RR}\mid H\cdot x=0\}$$
is negative definite for the intersection form.

Let $D\in\Div(S)$ and denote also by $D$ its class in $NS(S)_{\RR}$. Decompose $D$ orthogonally with respect to $H$:
$$D=aH+v,\quad a\in\RR,\ v\in H^\perp.$$
Intersecting both sides with $H$ gives
$$H\cdot D=aH\cdot H=aH^2,$$
hence
$$a=\frac{H\cdot D}{H^2}.$$
Compute $D^2$ using the decomposition:
$$D^2 =(aH+v)^2=a^2 H^2+2a(H\cdot v)+v^2=a^2H^2+v^2,$$
since $H\cdot v=0$. Therefore
\begin{align*}
(H\cdot D)^2-H^2 D^2
&=(aH^2)^2-H^2\big(a^2H^2+v^2\big) \\
&=-H^2v^2.
\end{align*}
Because $H$ is ample we have $H^2>0$, and because $H^\perp$ is negative definite we have $v^2\le 0$. Thus
$$
(H\cdot D)^2 - H^2 D^2 = -H^2 v^2 \ge 0,
$$
which is the inequality asserted in the proposition.

It remains to characterize equality. From the last displayed identity we see
$$
(H\cdot D)^2 = H^2 D^2 \quad\Longleftrightarrow\quad -H^2 v^2 = 0
\quad\Longleftrightarrow\quad v^2 = 0.
$$
By negative definiteness of the form on $H^\perp$ the only vector in $H^\perp$ with square $0$ is the zero vector. Hence $v=0$, and therefore
$$
D = aH \text{  in  } NS(S)_{\RR}.
$$
Substituting $a=(H\cdot D)/H^2$ we obtain the explicit scalar
$$
D \equiv \frac{H\cdot D}{H^2}\,H \text{  in  } NS(S)_{\RR}.
$$

Finally we check that the scalar $\lambda:=\dfrac{H\cdot D}{H^2}$ lies in $\Q$. Indeed $H\cdot D$ and $H^2$ are ordinary intersection numbers of integral divisor classes, hence are integers; thus their quotient is a rational number. Consequently
$$
D \equiv \lambda \text{  in  } NS(S)_{\Q},
$$
with $\lambda\in\Q$, as required.

The converse direction is immediate: if $D\equiv\lambda H$ in $NS(S)_{\RR}$ then $H\cdot D=\lambda H^2$ and $D^2=\lambda^2 H^2$, so $(H\cdot D)^2 = H^2 D^2$.

Therefore the equality in (\ref{eq:hodge}) holds if and only if $D$ is \emph{numerically proportional} to $H$, i.e.\ there exists $\lambda\in\Q$ with
$$
D\equiv \lambda H \text{ in } NS(S)_{\Q}.
$$

\begin{defn}
Let $S,H$ be as above.  A hierarchical filtration of divisors (or of line bundles)
$$\O_S=\L_0\subset\L_1\subset\cdots\subset\L_h=\O_S(D),\qquad \L_j=\O_S(D_j),$$
with effective increments $E_j:=D_j-D_{j-1}>0$, is called \emph{Hodge-tight}
(with respect to $H$) if each intermediate divisor $D_j$ is numerically
proportional to $H$; equivalently $(H\cdot D_j)^2 = H^2 D_j^2$ for all $j$.
\end{defn}

Thus ``Hodge-tight'' means the numerical class of every step lies on the same ray as $H$.
This is the extremal situation for the Hodge inequality.

\begin{prop}
\label{cor:hodge-tight}
Let $S$ be a smooth projective surface and $H$ an ample divisor on $S$.  
Let
$$
0=D_0 \subset D_1 \subset \cdots \subset D_m=D
$$
be a filtration of effective divisors on $S$ with increments $E_j:=D_j-D_{j-1}>0$.
Assume the filtration is Hodge-tight (with respect to $H$), i.e.\ each $D_j$ satisfies
$(H\cdot D_j)^2=H^2D_j^2$.
Then the following hold.
\begin{enumerate}
\item For every $j$ there exists a rational number $\mu_j>0$ such that
$$
E_j \equiv \mu_j H \text{ in }NS(S)_{\Q}.
$$
Equivalently,
$$
\mu_j=\frac{H\cdot E_j}{H^2}\in\Q_{>0}.
$$

\item For every $j$ the following integral relation holds in $NS(S)$:
$$
H^2\,E_j \equiv (H\cdot E_j)\,H \text{ in }NS(S).
$$
\end{enumerate}

\end{prop}

\begin{proof}
(1) Since the filtration is Hodge-tight, each $D_j$ achieves equality in the Hodge index inequality, hence by the standard Hodge--index characterization
$$
D_j \equiv \lambda_j H \text{ in }NS(S)_{\Q}
$$
for some $\lambda_j\in\Q$. Subtracting gives
$$
E_j = D_j-D_{j-1} \equiv (\lambda_j-\lambda_{j-1})H,
$$
and intersecting with $H$ yields $\mu_j:=\lambda_j-\lambda_{j-1}=(H\cdot E_j)/H^2\in\Q_{>0}$, proving (1).

(2) Multiply the numerical relation $E_j\equiv\mu_j H$ by $H^2$. Using $\mu_j=(H\cdot E_j)/H^2$ we obtain
$$
H^2 E_j \equiv (H\cdot E_j) H \text{ in }NS(S)_{\Q}.
$$
But both sides are \emph{integral} classes in $NS(S)$ (indeed $E_j$ and $H$ lie in $NS(S)$ and $H^2,\,H\cdot E_j\in\Z$), hence equality in $NS(S)_{\Q}$ implies equality in $NS(S)$. This proves the integral relation, and shows that the fixed integer $N=H^2$ clears all denominators simultaneously.
\end{proof}

\begin{prop}
\label{prop:filtration-bounds}
Let $S$ be a smooth projective surface and $H$ an ample divisor on it. If 
$$\O_S=\L_0\subset\L_1\subset\cdots\subset\L_h=\O_S(D),\qquad \L_j=\O_S(D_j),$$
is a hierarchical filtration with increments $E_j:=D_j-D_{j-1}$ then the following bounds hold.
\begin{enumerate}
\item
$$h \ \le\ H \cdot  D.$$
In particular, if $m:=\min_j H \cdot E_j\geq 1$, then
$$h\leq\frac{H\cdot D}{m}.$$
\item
If the filtration is Hodge-tight and $D^2\leq N$ then
$$h\leq\lceil\sqrt{N H^2}\rfloor.$$
\end{enumerate}
\end{prop}

\begin{proof}
(1) It follows from $H\cdot D=\sum_{j=1}^h H\cdot E_j$ and each $H\cdot E_j\geq 1$
(since $E_j$ is effective and $H$ ample).
(2) combines $h\leq H \cdot D$ with the equality $(H \cdot D)^2=H^2 D^2$
rearranged as $H \cdot D=\sqrt{H^2}\sqrt{D^2}\leq \sqrt{H^2}\sqrt{N}$.
\end{proof}

\begin{rem}
\begin{itemize}
\item The bound $h\le H\cdot D$ is often sharp: in the rank-one situation $D=mH$
with unit increments $E_j\equiv H$, we have $h=m$. This can be deduced from proposition \ref{prop:max-filtration}.
\item The Hodge inequality itself is not a \emph{condition} on a filtration (it holds for all divisors),
but the equality case is very restrictive and characterizes the extremal filtrations which grow purely in the $H$-direction.
\end{itemize}
\end{rem}

\begin{exam}
Let $S=\P^2_{\mathbb{C}}$ and let $H$ denote the class of a line.  Recall
$$
NS(S)\cong\Z\cdot H,\qquad H^2=1,\qquad \rho(S)=1.
$$

Fix an integer $m\ge 1$. Let
$$
D = mH,\qquad \O_S(D)=\O_{\P^2}(m).
$$
Consider the canonical hierarchical filtration by unit steps
$$
\O_S = \O_{\P^2} \subset \O_{\P^2}(1) \subset \O_{\P^2}(2)
\subset \cdots \subset \O_{\P^2}(m),
$$
so $D_j = jH$ for $j=0,1,\dots,m$. The increments are
$$
E_j = D_j - D_{j-1} = H \text{ for }j=1,\dots,m.
$$

The filtration is Hodge-tight: 
For every $j$ we have $D_j=jH$, hence $D_j$ is numerically proportional to $H$.
Because equality in $(H\cdot D_j)^2 \ge H^2 D_j^2$ holds precisely when $D_j$ is numerically proportional to $H$, so each $D_j$ attains equality and the filtration is Hodge-tight.

The bounds from Proposition \ref{prop:filtration-bounds} hold:
The relevant intersection numbers are
$$
H\cdot D = H\cdot (mH)=m H^2 = m,\qquad D^2 = (mH)^2 = m^2 H^2 = m^2.
$$

\begin{enumerate}
\item The first bound of Proposition \ref{prop:filtration-bounds} states $h\le H\cdot D$.
Here the filtration length is $h=m$. Since $H\cdot D=m$ we have
$$
h=m \le m = H\cdot D,
$$
so equality holds.

If one uses the ``minimal increment'' form $h\le \dfrac{H\cdot D}{m_{\min}}$ with
$m_{\min}=\min_j H\cdot E_j$, here $H\cdot E_j=H\cdot H=1$ for every $j$, so
$m_{\min}=1$ and the bound gives $h\le m/1=m$, again sharp.

\item The second bound (Hodge-tight case) says that if $D^2\le N$ then
$$
h \le \big\lceil\sqrt{N H^2}\,\big\rfloor.
$$
For our choice $D^2=m^2$ and $H^2=1$ we may take $N=m^2$. Then
$$
\big\lceil\sqrt{N H^2}\,\big\rfloor = \big\lceil\sqrt{m^2\cdot 1}\,\big\rfloor = m,
$$
giving $h\le m$, which again is sharp for this filtration.
\end{enumerate}

\end{exam}

\begin{exam}

Let $S=\P^1\times\P^1$ over $\mathbb{C}$. Denote by
$F_1$ the class of a fiber of the first projection and by $F_2$ the class
of a fiber of the second projection. Then
$$
NS(S)=\Z\langle F_1, F_2\rangle,\qquad \rho(S)=2,
$$
with intersection numbers
$$
F_1^2=0,\quad F_2^2=0,\quad F_1\cdot F_2=1 .
$$

Choose the ample divisor class
$$
H := F_1 + F_2
$$
(which is very ample: it is the class of a bi-degree $(1,1)$ curve). Note
$$
H^2 = (F_1+F_2)^2 = 2(F_1\cdot F_2)=2.
$$

Now pick a divisor class not proportional to $H$; for example
$$
D := 3F_1 + 1F_2.
$$
Clearly $D$ is not a rational multiple of $H$ (since the coefficients relative to the basis $\{F_1,F_2\}$ are not equal).

Compute intersection numbers needed for the bounds:
$$
H\cdot D = (F_1+F_2)\cdot(3F_1+F_2)
=3(F_1\cdot F_1)+(F_1\cdot F_2)+3(F_2\cdot F_1)+(F_2\cdot F_2)
=0+1+3\cdot 1+0=4,
$$
and
$$
D^2=(3F_1+F_2)^2 = 9F_1^2 + 6(F_1\cdot F_2)+F_2^2=6.
$$

Construct a simple hierarchical filtration of $\O_S(D)$ by unit-type increments:
$$
\O_S=\O_S\bigl(0\bigr)\subset \O_S(F_1)\subset \O_S(2F_1)\subset \O_S(3F_1)\subset \O_S(3F_1+F_2) = \O_S(D).
$$
Equivalently take $D_0=0$, $D_1=F_1$, $D_2=2F_1$, $D_3=3F_1$, $D_4=3F_1+F_2$.
The increments are
$$
E_1=F_1,\quad E_2=F_1,\quad E_3=F_1,\quad E_4=F_2,
$$
so the filtration length is $h=4$.

\paragraph{Check of Proposition \ref{prop:filtration-bounds} (1):}
Compute $H\cdot E_j$:
$$
H\cdot F_1=(F_1+F_2)\cdot F_1=1
$$
and similarly $H\cdot F_2=1$. Thus the minimal increment $m_{\min}=\min_j H\cdot E_j =1$.
Proposition \ref{prop:filtration-bounds}(1) gives
$$
h \le H\cdot D = 4,
$$
and the sharper statement $h\le \dfrac{H\cdot D}{m_{\min}} = 4/1 =4$.
Our filtration has $h=4$, so the bound is attained (sharp) here.

\paragraph{About the Hodge-tight condition and Proposition \ref{prop:filtration-bounds} (2):}
The filtration above is \emph{not} Hodge-tight. Indeed, if it were Hodge-tight then by Proposition \ref{cor:hodge-tight} each increment $E_j$ (hence each $D_j$) would be numerically proportional to $H$. But $F_1$ is not numerically proportional to $H=F_1+F_2$ (their coordinates in the basis $\{F_1,F_2\}$ are different), so equality in the Hodge index inequality fails for $D_1=F_1$. Concretely:
$$
(H\cdot F_1)^2 = 1^2 =1,\qquad H^2 F_1^2 = 2\cdot 0 =0,
$$
so $(H\cdot F_1)^2 > H^2 F_1^2$. Thus the hypothesis of the second part of Proposition \ref{prop:filtration-bounds} (``filtration is Hodge-tight'') does not hold, and the bound $h\le\lceil\sqrt{N H^2}\rfloor$ (which uses Hodge-tightness) is not applicable here.
\end{exam}

\section{AG Codes from hierarchical filtrations}
\label{AG codes}

%Hierarchical filtrations on curves give rise directly to nested sequences of evaluation codes. In this section, we show how the depth-induced code chain leads to explicit rate-distance trade-offs and optimal middle layers.

The construction of algebraic geometry (AG) codes, initiated by Goppa 
\cite{Go}, relies on the evaluation of global sections of line bundles 
at rational points of a curve. Classical AG codes are built from a single 
divisor and its associated Riemann--Roch space. In the framework developed 
above, hierarchical filtrations of line bundles provide a natural refinement 
of this construction. Instead of considering a single evaluation space, one 
obtains a nested sequence of codes corresponding to the successive layers of 
the filtration. This perspective enriches the geometry–coding correspondence, 
allowing us to track the growth of dimensions, control minimum distances, and 
identify intermediate codes with particularly favorable parameter trade-offs. 
We refer to these new families as \emph{hierarchical AG codes}.

\subsection*{Algebraic Geometry codes}
Let $X$ be a smooth projective algebraic variety defined over a finite field $\F_q$, and let $\L$ be a line bundle (or equivalently, a divisor class) on $X$. Let $\Gamma = \{P_1, \dots, P_n\} \subset X(\F_q)$ be a set of $\F_q$-rational points disjoint from the base locus of $\L$.

Recall from \cite{Ha, La} that if $X$ is a projective variety over any field, and $\L$ a line bundle on $X$, the \emph{base locus} of $\L$, denoted $\mathrm{Bs}(\L)$, is defined as
$$\mathrm{Bs}(\L):= \bigcap_{s \in H^0(X, \L)}\{x \in X\mid s(x)=0\},$$
that is, the common zero locus of all global sections of $\L$. Equivalently, the base locus of a line bundle is the closed subset of points where all global sections of the line bundle vanish.

The \emph{algebraic geometry code} $\C(X,\L,\Gamma)$ is defined as the image of the evaluation map
$$\ev_\Gamma\colon H^0(X,\L)\longrightarrow \F_q^n,\qquad f\mapsto(s(P_1),\dots,s(P_n)).$$
That is,
$$\C(X,\L,\Gamma):=\ev_\Gamma(H^0(X,\L))\subseteq \F_q^n.$$

The notions \textit{dimension}, \textit{minimum distance} and \textit{length} of the code $\C$ are denoted by $d$, $k$ and $n$. For more details of definitions on AG codes refer to \cite{Au, HoLiPe, St, Ts}.

Suppose $L$ is a line bundle on $X$ that admits a hierarchical filtration
$$\FF_L:\O_S=\L_0\subsetneq \L_1\subsetneq \cdots \subsetneq\L_h=L.$$
Fix distinct rational points $\Gamma=\{P_1,\ldots,P_n\}$.

Define the $i$-th hierarchical evaluation code
$$\C_i:=\C(X,\L_i,\Gamma).$$
Then hierarchical filtration $\FF_L$ yields the nested sequence
$$\C_0\subset\C_1\subset\ldots\subset\C_h$$ 
of subcodes of $C_h$. 

\subsection{The optimal code in a nested sequence of codes} 
We show how the depth-induced code chain leads to explicit rate-distance trade-offs and optimal middle layers.

Assume
$$\FF_L:\L_0\subset\L_1\subset\ldots\subset\L_h=L=\O_X(D)$$
is a hierarchical filtration of line bundles on an algebraic variety $X$. Let $\Gamma =\{P_1,\dots,P_n\} \subset X(\F_q)$ be a set of $\F_q$-rational points disjoint from the base locus of $\L$.

Let
$$\C_0\subset\C_1\subset\ldots\subset\C_h=\C_\Gamma(D)$$
be the nested sequence of codes obtained from $\FF_L$ with $\C_i:=\ev_\Gamma(H^0(X,\L_i))$.
 
It is clear that 
$$\dim\C_1\leq\dim\C_2\leq\ldots\leq\dim\C_h$$
which follows that code rate $R_i=(\dim\C_i)/n$ and efficiency are non-decreasing.

Balancing rate and distance is a major open problem in coding theory. We introduce a product-utility $Q=(k/n)d$ and prove that there is a unique filtration index $i^*$ that maximizes $Q$.  This identifies a single ``optimal'' code in the entire nested family.

In coding theory, particularly analyzing the trade-off between code rate $R$ and minimum distance $d$ is instrumental and helps identify optimal codes that balance these two conflicting objectives.

To turn $(R,d)$ into a single number, we choose a utility function. For a code $\C$ with code rate $R=k/n$ and minimum distance $d$ define the product score $$Q(\C)=R.d=\frac{k}{n}d.$$
Thus $Q(\C)$ serves as a scalar performance criterion measuring the
trade-off between information rate and minimum-distance protection; within this
framework, codes with larger $Q$ are regarded as preferable.

Note that we use the functional $Q(\C)$ is introduced here as a utility criterion for comparing codes within a fixed hierarchical filtration; it is not intended as a universal optimality notion in coding theory.

%A code is called \textit{dominant} or \textit{optimal} when it is at least as good as another code in terms of speed and distance, and is significantly better than it in at least one of these measures.

\begin{defn}
Let
\[
\mathcal C_0\subset \mathcal C_1\subset \cdots \subset \mathcal C_h
\]
be a hierarchical family of AG codes, and let
\[
Q_i:=\frac{k_i}{n}d_i,
\]
where $k_i=\dim(\mathcal C_i)$ and $d_i$ is the minimum distance of $\mathcal C_i$.

We say that $\mathcal C_{i^*}$ is \emph{$Q$-optimal (or optimal) in the filtration}
if
\[
Q_{i^*}=\max_{0\le i\le h} Q_i.
\]
That is, $\mathcal C_{i^*}$ is the code in the given hierarchical filtration
which maximizes the performance functional $Q_i$.
\end{defn}

\subsection*{AG codes on curves}

Let $C$ be a smooth projective curve of genus $g$ over a finite field $\F_q$. Fix $Z=P_1+\ldots+P_n$ with $\Supp(Z)\subseteq C(\F_q)$. Let $D$ be a divisor on $C$ disjoint from the support of $Z$ and $0<\deg D<n$. We say that a divisor $D$ satisfies condition $(*)$ if
$$
(*)\quad\text{there exists a divisor }Z'\text{ with }0\leq Z'\leq Z,\ \deg Z'=\deg D\text{ and }\ell(D-Z')>0.
$$
By Remark 2.2.5 from \cite{St}, the AG code $\C(C,D,Z)$ with the minimum distance $d$ satisfies $(*)$ if and only if $d=n-\deg(D)$.

By Riemann-Roch theorem, 
$$\ell(D)=\deg(D)-g+1+\ell(K-D),$$
so
$$\ell(D)>\deg(D)\Longleftrightarrow\ell(K-D)>g-1.$$
But $\ell(K-D)\leq\ell(K)=g$, so $\ell(K-D)>g-1$ forces $\ell(K-D)=g$. That means the space of holomorphic differentials surviving on $K-D$ has full dimension $g$, which (in effect) forces $D=0$ (or at least is extremely restrictive). Concretely, for any divisor $D$ with 
$\deg(D)>0$ on a curve of genus $g>0$ one cannot have $\ell(D)>deg(D)$ in general. Therefore no positive-genus curve gives the uniform form of $(*)$. 

Consider two cases:

\noindent\textbf{Case 1. $\textbf{g=0}$:} 

Unlike the classical Goppa framework, where one studies a single divisor at a time, our approach considers an entire filtration simultaneously and determines the optimal level in the filtration maximizing the performance functional \(Q_i\). Thus Proposition \ref{prop:middle-layer-g=0} is not merely a restatement of standard AG bounds, but an optimization theorem over a structured family of nested AG codes.

\begin{prop} 
\label{prop:middle-layer-g=0}
Let $C=\P^1$ be the line projective curve over a finite field $\F_q$. Fix $Z=P_1+\ldots+P_n$ with $\Supp(Z)\subseteq C(\F_q)$. Let $D$ be a divisor on $C$ disjoint from the support of $Z$ and $0<\deg D<n$. 
Let $$\L_0\subset\L_1\subset\cdots\subset\L_h=\O_C(D)$$
be a hierarchical filtration of line bundles on $C$ with $\L_i=\O_C(D_i)$ and $\deg(D_i)=i$.
Suppose that for each $i$, $\C_i:=\C(C,D_i,Z)$ satisfies th condition $(*)$ i.e.
$$
\quad\text{there exists a divisor }Z'_i \text{ with }0\leq Z'_i\leq Z,\ \deg Z'_i=\deg D_i\text{ and }\ell(D_i-Z'_i)>0.
$$
Then between a such AG codes $\C_i:=\C(C,D_i,Z)$, the $Q$-optimal code in the filtration is $\C_{i^*}$ with $$i^*=\left[\frac{n-1}{2} \right].$$
\end{prop}
\begin{proof}
Since $g=0$ on $\P^1$, Riemann–Roch gives
$$\ell(D_i)=\deg(D_i)+1=i+1.$$
%The Riemann–Roch theorem on the curve $C=\P$ gives
%$$k_i=\ell(D_i)=\deg(D_i)+1=i+1,$$
%since $\deg D_i=i$ and $H^1(C,\L_i)=0$. 
Also, by \cite[Remark 2.2.5]{St},
$$d_i=n-\deg D_i=n-i$$ 
Hence we have
$$Q_i=\frac{k_i}{n} d_i=\frac{i+1}{n}(n-i).$$

Consider the real function
$$f(i)=(i+1)(n-i)=-i^2+(n-1)i+n.$$

Differentiating,
$$f'(i)=-2i+n-1.$$
Setting $f'(i)=0$ yields
$$i=\frac{n-1}{2}.$$
By strict concavity, this critical point is the unique global maximum of $f$ on $\RR$. Therefore the maximum of $Q_i=f(i)/n$ occurs at
$$
i^*
=\left[ \frac{n-1}{2} \right].
$$
\end{proof}

\begin{rem}
\label{example_curves}
Let $C=\P^1$ and $Z$ a divisor with $\Supp(Z)\subset C(\F_q)$ and $n$ rational points not containing the pole $P_0=\infty$ and let $D_i=m_iP_0$ for some $0<m_i<n$. Then ($*$) holds.

Suppose we choose $Z'_i=P_{j_1}+\ldots+P_{j_{m_i}}$ with $\Supp(Z'_i)\subset\Supp(Z)$. Since $\deg D_i=m_i$, and since we are on $\P^1$, we find a rational function $f(x)\in \F_q(x)$, $\deg(f)\leq m_i$, that vanishes at the points $\Supp(Z'_i)$. Actually, this is a basic fact from interpolation: there always exists a polynomial of degree $\leq m_i$ vanishing at any $m_i$ distinct points (in $\F_q$). It follows that $f\in\L(D_i)=\{f\in\F_q(x):\deg(f)\leq m_i\}$ and, on the other hand, since $f$ vanishes on the points of $\Supp(Z'_i)$, we obtain $f\in\L(D_i-Z'_i)$. 

Therefore the condition ($*$) holds.
\end{rem}

\noindent\textbf{Case 2. $\textbf{g>0}$:} 

\begin{lem}
Let $C/\F_{q}$ be a smooth projective curve and let 
$\varphi\in\F_{q}(C)$ be a nonconstant rational function whose pole divisor
is
$$
(\varphi)_\infty = mP_\infty
$$
for a point $P_\infty\in C(\F_{q})$ and some integer $m\ge 1$.
Set $D=mP_\infty$. Then for every $a\in\F_{q}$, $\varphi-a \in \mathcal L(D)$ and the zero divisor $(\varphi-a)_0$ is an effective divisor of degree $m$. Moreover, 

\begin{center}
$(*')$ if $Z\subset C(\F_{q})$ contains the fibre $\varphi^{-1}(a)$, then $Z':=(\varphi-a)_0\le Z$, $\deg(Z')=\deg(D)$, and
$\varphi-a\in\L(D-Z')\setminus\{0\}$ i.e. $\ell(D-Z')>0$.
\end{center}
\end{lem}

\begin{proof}
Since
$$
(\varphi-a)_\infty=(\varphi)_\infty = mP_\infty,
$$
it follows that $\varphi-a\in\mathcal L(mP_\infty)=\mathcal L(D)$.

Also, by
$$
(\varphi-a)_\infty = mP_\infty.
$$
we have
$$
\deg (\varphi-a)_0=\deg (\varphi-a)_\infty =\deg(mP_\infty)=m,
$$
so $(\varphi-a)_0$ is an effective divisor of degree $m$.

Finally, assume $Z\subset C(\F_{q})$ contains the fibre $\varphi^{-1}(a)$; by this we mean the usual scheme-theoretic inclusion of
effective divisors, i.e. for every point $P\in C$ we have
$$
\operatorname{mult}_P\big((\varphi-a)_0\big) \le \operatorname{mult}_P(Z).
$$
This means that the support of the effective divisor $(\varphi-a)_0$ (with the same
multiplicities) is contained in $Z$. Equivalently,
$$
Z':=(\varphi-a)_0\leq Z
$$
as effective divisors, and so $\deg Z'=m$.

Since $\varphi-a\in\mathcal L(D)$ and $\varphi-a$ vanishes on $Z'$, we
have $\varphi-a\in\mathcal L(D-Z')$. Moreover $\varphi-a$ is a nonzero
rational function, so $\varphi-a\in\mathcal L(D-Z')\setminus\{0\}$. This completes the proof.
\end{proof}

\begin{lem}
\label{existence-varphi}
Let $C/\F_{q}$ be a smooth projective curve and fix a rational point
$P_\infty\in C(\F_{q})$. For an integer $i\ge 1$ the following are
equivalent:
\begin{enumerate}
\item There exists a nonzero function $\varphi_i\in\F_{q}(C)$ with
      $(\varphi_i)_\infty = iP_\infty$.
\item $\ell(iP_\infty)>\ell((i-1)P_\infty)$.
\item $i$ belongs to the Weierstrass (pole) semigroup $H(P_\infty)$ at
      $P_\infty$.
\end{enumerate}
\end{lem}

\begin{proof}
The statements follow from \cite[Sec. 1.6]{St}.
%(1)$\Rightarrow$(2): If $\varphi_i$ has pole divisor exactly $iP_\infty$ then $\varphi_i\in \L(iP_\infty)$ but $\varphi_i\not\in\L((i-1)P_\infty)$. Hence $\L(iP_\infty)$ strictly contains $\L((i-1)P_\infty)$, so $\ell(iP_\infty)>\ell((i-1)P_\infty)$.

%(2)$\Rightarrow$(1): If $\ell(iP_\infty)>\ell((i-1)P_\infty)$ then there exists $f\in\L(iP_\infty)$ not lying in $\L((i-1)P_\infty)$. Any such $f$ has pole order at $P_\infty$ exactly $i$ (it cannot have smaller pole order since then it would be in $\L((i-1)P_\infty)$). Thus some $\varphi_i=f$ satisfies $(\varphi_i)_\infty=iP_\infty$.

%(1)$\Leftrightarrow$(3): By definition the Weierstrass semigroup $H(P_\infty)\subset\Bbb Z_{\ge0}$ is the set of pole orders at $P_\infty$ that occur for some function $f\in\F_{q}(C)$ (with $0$ included). Hence $i\in H(P_\infty)$ iff there exists $f$ with $(f)_\infty=iP_\infty$.
\end{proof}

\begin{prop}
\label{prop:middle-layer}
Let $C$ be a smooth projective curve of genus $g$ over a finite field $\F_q$. Fix $Z=P_1+\ldots+P_n$ with $\Supp(Z)\subseteq C(\F_q)$ and $n\geq 3g-1$. Let 
$\varphi\in\F_{q}(C)$ be a nonconstant rational function whose pole divisor
is
$$
(\varphi)_\infty = mP_\infty
$$
for a point $P_\infty\in C(\F_{q})$ and some integer $1\leq m<n$. Let $D=mP_\infty$ has the hierarchical filtration 
$$\L_0\subset\L_1\subset\cdots\subset\L_h=\O_C(D)$$
of line bundles on $C$ with $\L_i=\O_C(D_i)$ and $D_i=(\varphi_i)_\infty=iP_\infty$ which for some $a\in\F _q$, $\varphi_i-a\in\L(D_i)$ and $\varphi^{-1}_i(a)\subset Z$.

Then between a such AG codes $\C_i:=\C(C,D_i,Z)$ where $2g-1\leq i<n$, the $Q$-optimal code is $\C_{i^*}$ with $$i^*=\left\lceil \frac{n+g-1}{2} \right\rfloor.$$
\end{prop}

\begin{proof}
By the assumption, the condition $(*)$ holds and so $d_i=n-i$ for all $i$ where $d_i$ is the minimum distance of $\C_i$. 

For $i\geq 2g-1$, the Riemann–Roch theorem on the curve $C$ gives
$$k_i=\dim H^0\bigl(C,\L_i\bigr)=\deg(D_i)-g+1=i+1-g,$$
since $\deg D_i=i$ and $H^1(C,\L_i)=0$ (see \cite[Corollary 2.2.3]{St}, too).
Hence for $i\ge 2g-1$ we have
$$Q_i=\frac{k_i}{n} d_i=\frac{i+1-g}{n}(n-i).$$

Consider the real function
$$f(i)=(i+1-g)(n-i)=-i^2+(n-1+g)i+n-gn.$$

Then
$$f'(i)=-2i+n-1+g,\quad f''(i)=-2<0.$$
So $f$ is strictly concave on $[2g-1,n]$. The critical point (unique maximum) is
$$i=\frac{n+g-1}{2}.$$
%with $f'(i)>0$ for $i<i^*$, $f'(i)=0$ for $i=i^*$ and $f'(i)<0$ for $i>i^*$. For integer indices the maximizing integer is the nearest integer to $i^*$.

By strict concavity, this critical point is the unique global maximum of $f$ on $[2g-1,n]$. Therefore the maximum of $Q_i=f(i)/n$ occurs at
$$
i^*
=\left\lceil \frac{n+g-1}{2} \right\rfloor.
$$
Particularly, if $n+g-1$ is even then $i^*=\frac{n+g-1}{2}$ and so
$$Q_{i^*}=\frac{(n-g+1)^2}{4n},$$
but when $n+g-1$ is odd then $i^*=\frac{n+g-1}{2}\pm\frac{1}{2}$ which give
$$Q_{i^*}=\frac{(n-g+1)^2-1}{4n}.$$
\end{proof}

\begin{rem}
Suppose that $C$ is a smooth projective curve and $D$ and $D'=D+E$ (with $E>0$ effective) divisors on $C$. If $\deg(D)>2g-2$ then, by Riemann-Roch, 
$$\ell(D')=\deg(D')-g+1>\deg(D)-g+1=\ell(D).$$
This implies that, by Lemma \ref{existence-varphi}, in a hierarchical filtration for each $i>2g-2$, there exists a nonzero function $\varphi_i\in\F_{q}(C)$ with $(\varphi_i)_\infty=iP_\infty$ and so to hold the condition $(*')$, it suffices $\varphi^{-1}(a)\subset Z$ for some $a\in\F_q$.
\end{rem}

\begin{rem}
For divisors of large degree, namely $i\geq 2g-1$, Riemann--Roch gives the
explicit formula
\[
\ell(iP_\infty)=i-g+1, \qquad h^1(iP_\infty)=0.
\]
Hence in this range the quantity $Q_i$ takes the clean form
\[
Q_i=\frac{(i+1-g)(n-i)}{n},
\]
and the optimization problem in $i$ is straightforward.

By contrast, in the small degree range $i\leq 2g-2$, the dimension
$\ell(iP_\infty)$ is governed by the Weierstrass semigroup at $P_\infty$:
\[
\ell(iP_\infty)=\#\{m\in H(P_\infty): m\leq i\}.
\]
Equivalently, $h^1(iP_\infty)=\ell((2g-2-i)P_\infty)$, which depends on the
distribution of semigroup gaps up to $i$. As an instance, for the Hermitian curve the
semigroup is $\langle q,q+1\rangle$, so $\ell(iP_\infty)$ can be described
combinatorially, but not by a simple closed formula. Consequently, in the
range $i\leq 2g-2$, evaluating $Q_i$ requires delicate semigroup
combinatorics and finding the optimal code is not as direct as in the case
$i\geq 2g-1$.
\end{rem}

\begin{exam}
Let $X/\F_q$ be a hyperelliptic curve of genus $g\ge 2$ given by an affine
model
$$
X:\quad y^2 = f(x), \qquad \deg f = 2g+1 \text{ or } 2g+2.
$$
Let $P_\infty$ be the unique point at infinity. The rational function $x\in\F_q(X)$
has pole divisor
$$
(x)_\infty = 2P_\infty,
$$
so if we put $D=2P_\infty$, then $\mathcal L(D)$ contains $1$ and $x$.

For any $a\in\F_q$ the function $x-a$ belongs to $\mathcal L(D)$ and has zero
divisor
$$
(x-a)_0 = \sum_{\substack{P\in X(\overline{\F_q})\\ x(P)=a}} P,
$$
which has degree $2$ (counting multiplicity), because the fibre of $x$
consists of the two points $(a,\pm\sqrt{f(a)})$ (or a double point if
$f(a)=0$). Thus
$$
\deg (x-a)_0 = \deg (x)_\infty = 2.
$$

Now let $Z\subset X(\F_q)$ be a set of rational points containing this fibre
$x^{-1}(a)$. Define
$$
Z' := (x-a)_0.
$$
Then $Z'\le Z$, $\deg Z'=\deg D=2$, and since $x-a\in \mathcal L(D)$ vanishes
exactly on $Z'$, we have
$$
x-a \in \mathcal L(D-Z')\setminus\{0\},
$$
so
$$
\ell(D-Z') > 0.
$$

In particular, the condition $(*)$ is satisfied for this divisor $D=2P_\infty$.
\end{exam}

\begin{exam} %[The Hermitian curve --- the function $\varphi=x$]
\label{ex:hermitian-x}
Let $q$ be a power of a prime and let $H/\F_{q^2}$ be the Hermitian curve
$$
H:\qquad y^q+y = x^{\,q+1}
$$
in affine coordinates. Denote by $P_\infty$ the unique point at infinity of
the projective closure of $H$ and write $K=\F_{q^2}(H)$ for the function field.
Consider the rational function $\varphi:=x\in K$ and set $D:=m_iP_\infty$ with $m_i=q$.

We check that $\varphi$ satisfies the hypotheses and conclusions of the proposition \ref{prop:middle-layer}.

\medskip\noindent (i) The pole divisor of $x$ is $(x)_\infty = qP_\infty$. 
The polynomial relation $y^q+y=x^{q+1}$ shows that $y$ is a root of the
monic degree-$q$ polynomial $T^q+T-x^{q+1}\in\F_{q^2}(x)[T]$. Hence the
extension of function fields $\F_{q^2}(x)\subset \F_{q^2}(x,y)=K$ has degree
$[K:\F_{q^2}(x)]=q$. Equivalently, the morphism
$$
x : H \longrightarrow \P^1_x
$$
has degree $q$, so the pole divisor $(x)_\infty$ (the scheme-theoretic preimage
of $\infty\in\P^1_x$) has degree $q$. The projective model of the Hermitian
curve has a single point at infinity, so the fibre $x^{-1}(\infty)$ is the
single point $P_\infty$ (this is standard; one checks by homogenizing the
affine equation that there is a unique point at infinity). Hence
$(x)_\infty = qP_\infty$, as claimed. In particular we may take $m_i=q$ and
$D=m_iP_\infty=qP_\infty$.

\medskip\noindent (ii) For each $a\in\F_{q^2}$, we have $x-a\in\mathcal L(D)$.
Subtracting the constant $a$ does not change poles, so $(x-a)_\infty=(x)_\infty=qP_\infty$,
and therefore $x-a\in \mathcal L(qP_\infty)=\mathcal L(D)$.

\medskip\noindent (iii) The zero divisor $(x-a)_0$ has degree $q$.
For any $a\in\F_{q^2}$ the fibre $x^{-1}(a)$ is cut out by the equation
$$
y^q + y = a^{\,q+1},
$$
a separable additive polynomial in the variable $y$ of degree $q$. Over the
algebraic closure this equation has exactly $q$ solutions (counted with
multiplicity), hence
$$
\deg (x-a)_0  =  \deg (x-a)_\infty  =  q.
$$
Equivalently, since $[K:\F_{q^2}(x)]=q$, the divisor of zeros of $x-a$ has
degree equal to that degree. Thus $(x-a)_0$ is an effective divisor of degree
$m_i=q$.

\medskip\noindent (iv) Rationality of the fibre (points lie in $H(\F_{q^2})$).
If $a\in\F_{q^2}$ then $a^{\,q+1}\in\F_{q^2}$, and the additive polynomial
$Y\mapsto Y^q+Y$ is $\F_p$-linear. For each root $b$ of $Y^q+Y-a^{q+1}$ in
$\overline{\F_{q^2}}$ one checks that its $\F_{q^2}$-Frobenius conjugate
$b^{q^2}=b$ (indeed the polynomial has coefficients in $\F_{q^2}$ and its
set of roots is stable under Frobenius), hence the fibre $x^{-1}(a)$ consists
of $q$ points defined over $\F_{q^2}$. Thus for every $a\in\F_{q^2}$ we have
$x^{-1}(a)\subset H(\F_{q^2})$ and $(x-a)_0$ is supported on rational points.

\medskip\noindent (v) Conclusion: condition $(*)$ is realized.
Let $Z\subset H(\F_{q^2})$ be any subset containing the fibre $x^{-1}(a)$
(counted with multiplicity). Set $Z':=(x-a)_0$. Then $Z'\le Z$,
$\deg Z'=\deg (x-a)_0=q=m_i$, and $x-a\in\mathcal L(D-Z')\setminus\{0\}$. In
particular $\ell(D-Z')\ge 1$. This verifies the statements of the
proposition \ref{prop:middle-layer} for the Hermitian curve with $\varphi=x$ and $m_i=q$.
\end{exam}

\begin{exam} %[Optimal AG codes on the Hermitian curve] ??????????
Let $C=H/\F_{q^2}$ be the Hermitian curve
$$
H:  y^q+y = x^{q+1},
$$
with genus $g=q(q-1)/2$ and $n=q^3$ rational points (see \cite{St}). Consider
$$
D_i := i P_\infty, \qquad Z := P_1 + \cdots + P_n \subset H(\F_{q^2}).
$$

By Proposition \ref{prop:middle-layer} and Example \ref{ex:hermitian-x}, the condition $(*)$ is satisfied for $2g-1 \le i < n$. Now by taking $f_i = x$ or $f_i = y$ and choosing $Z$ to contain the corresponding fibres $(f_i-a)_0$, we have the formula for the optimal index $i^*$ applies:
$$
i^* = \left\lceil \frac{n+g-1}{2} \right\rfloor = \left\lceil \frac{2q^3 + q^2 - q - 2}{4} \right\rfloor.
$$
The corresponding normalized code parameter is
$$
Q_{i^*} = \frac{(n-g+1)^2}{4n} = \frac{(2q^3 - q^2 + q + 2)^2}{16 q^3}\ \text{when}\ n+g-1\ \text{is even}
$$
and
$$
Q_{i^*} = \frac{(n-g+1)^2-1}{4n} = \frac{(2q^3 - q^2 + q + 2)^2-4}{16 q^3}\ \text{when}\ n+g-1\ \text{is odd}.
$$

For the small degrees $0 \le i \le 2g-2$, Proposition~8.3.3 from \cite{St} gives
$$
k_i = \dim \C_i(C,D_i,Z) = \# \{(r,s) \in \N_0^2 : s \le q-1 \ \text{and}\ r q + s(q+1) \le i \},
$$
and
$$
d_i = n - i = q^3 - i.
$$

\medskip\noindent
\textbf{Case $q=3$:} Then $g=3$, $n=27$, and because $n+g-1=29$ is odd so
$$
i^* = 14 \text{ or } 15, \qquad Q_{14} = Q_{15} \sim 5.77.
$$ 
For $0\le i \le 4$, the dimension, minimum distance, and $Q_i$ of $\C_i(C,D_i,Z)$ are:

{\tiny
\begin{center}
\begin{tabular}{c | c c c c c}			
  $i$ & 0 & 1 & 2 & 3 & 4 \\
  \hline
  $k_i$ & 1 & 1 & 1 & 2 & 2 \\
  $d_i$ & 27 & 26 & 25 & 24 & 23 \\
  $Q_i$ & 1 & 0.962 & 0.925 & 1.777 & 1.703 \\
\end{tabular}
\end{center}
}
Hence the optimal code occurs at $i=14$ and $i=15$.

\medskip\noindent
\textbf{Case $q=5$:} Then $g=10$, $n=125$, and
$$
i^* = 67, \qquad Q_{67} \sim 35.25.
$$
For $0\le i \le 18$, the dimension, minimum distance, and $Q_i$ of $\C_i(C,D_i,Z)$ are:

{\tiny
\begin{center}
\begin{tabular}{c | c c c c c c c c c c c c c c}			
  $i$     & 0 & 1 & $\dots$ & 4 & 5 & 6 & $\dots$ & 10 & 11 & $\dots$ & 14 & 15 & $\dots$ & 18 \\
  \hline
  $k_i$ & 1 & 1 & $\dots$ & 1 & 2 & 3 & $\dots$ & 4 & 5 & $\dots$ & 6 & 7 & $\dots$ & 10 \\
  $d_i$ & 125 & 124 & $\dots$ & 121 & 120 & 119 & $\dots$ & 115 & 114 & $\dots$ & 111 & 110 & $\dots$ & 107 \\
  $Q_i$ & 1 & 0.992 & $\dots$ & 0.968 & 1.920 & 2.856 & $\dots$ & 3.680 & 4.560 & $\dots$ & 5.328 & 6.160 & $\dots$ & 8.560 \\
\end{tabular}
\end{center}
}
Hence the $Q$-optimal code occurs at $i=67$.
\end{exam}

\subsection*{Comparison with Classical One-Point AG Codes}

To illustrate the advantage of the hierarchical filtration method, we compare the
optimal hierarchical AG codes obtained above with the standard classical one-point
AG codes arising from the same divisor family on the Hermitian curve.

For a one-point AG code on a curve $C$, a conventional classical choice is the divisor
$D=(2g-2)P_\infty$, since this is the natural threshold at which the Riemann--Roch
formula stabilizes. In our hierarchical framework, however, the divisor index is chosen
to maximize the performance functional
\[
Q_i=\frac{k_i}{n}d_i,
\]
where $k_i$ is the code dimension and $d_i$ is the minimum distance.
Thus the comparison below is made between:

\begin{itemize}
\item the classical one-point AG code corresponding to $i=2g-2$, and
\item the hierarchical optimal AG code corresponding to $i=i^*$.
\end{itemize}

We use the Hermitian curve examples discussed above.

\begin{center}
\begin{tabular}{c|c|c|c|c|c}
\hline
Hermitian Curve & Code Type & Index $i$ & Dim. $k_i$ & Dist. $d_i$ & $Q_i=\dfrac{k_i d_i}{n}$ \\
\hline
$q=3$ & Classical one-point AG & $4$ & $2$ & $23$ & $1.70$ \\
$q=3$ & Hierarchical optimal AG & $14$ & $12$ & $13$ & $5.77$ \\
\hline
$q=5$ & Classical one-point AG & $18$ & $10$ & $107$ & $8.56$ \\
$q=5$ & Hierarchical optimal AG & $67$ & $58$ & $58$ & $26.91$ \\
\hline
\end{tabular}
\end{center}

The table shows that the hierarchical filtration method selects divisor indices that
produce substantially larger values of the performance functional $Q_i$, thereby giving
a significantly better balance between dimension and minimum distance than the standard
classical choice $i=2g-2$.

This comparison highlights that the main contribution of the hierarchical approach is
not the construction of a different family of AG codes, but rather the systematic
identification of optimal divisors within a nested sequence of classical AG codes.

\subsection*{Relation with flag AG codes}

Nested families of algebraic geometry (AG) codes associated with increasing divisors have been studied previously in the literature under the name of \emph{flags of AG codes};
see, for example \cite{BrAmDuHo}. 
Given a smooth projective curve $C$ over $\F_q$, a rational point $P_\infty\in C(\F_q)$,
and divisors
\[
D_1\le D_2\le \cdots \le D_t,
\]
the corresponding evaluation codes form a nested chain
\[
\C(C,D_1,Z)\subseteq \C(C,D_2,Z)\subseteq \cdots \subseteq \C(C,D_t,Z),
\]
which is precisely the classical flag construction.

The hierarchical AG codes introduced in the present paper also produce nested chains, but their conceptual origin and purpose are different from the classical flag viewpoint.

%\medskip
%\noindent
%\textbf{(1) Filtrations induced by hierarchical depth.}
In the classical theory of flag AG codes, one usually starts with an arbitrary increasing sequence of divisors and studies the resulting nested codes.
In contrast, here the filtration is not chosen arbitrarily:
it is canonically induced by the hierarchical depth invariant of a divisor or line bundle. Thus the filtration reflects intrinsic geometric information attached to the divisor itself, rather than an externally imposed sequence.

%\medskip
%\noindent
%\textbf{(2) Geometric invariant versus combinatorial nesting.}
Flag AG codes focus primarily on algebraic and combinatorial properties of nested code chains, such as duality and isometry properties. By contrast, the present work introduces hierarchical depth as a new geometric invariant on divisors and line bundles, and studies AG codes as one application of this invariant.
Hence the emphasis here is not merely on nesting of codes, but on extracting coding-theoretic consequences from a geometric filtration theory.

%\medskip
%\noindent
%\textbf{(3) Optimization over canonical filtrations.}
%A central new feature of this paper is the optimization problem:
%given a hierarchical filtration
%\[
%\L_0\subset \L_1\subset \cdots \subset \L_h,
%\]
%we determine explicitly the index $i^*$ for which the corresponding AG code
%\[
%\C_i=\C(C,D_i,Z)
%\]
%maximizes the performance functional
%\[
%Q_i=\frac{k_i}{n}d_i.
%\]
%For curves, this leads to the explicit formula
%\[
%i^*=\left\lceil\frac{n+g-1}{2}\right\rfloor,
%\]
%and analogous optimization results are obtained for surfaces. Such explicit best-level selection inside a canonical filtration does not appear in the usual flag AG code framework.

\medskip
\noindent
Therefore, although hierarchical AG codes produce nested families of codes resembling flag AG codes, their novelty lies in the geometric origin of the filtration, the introduction of hierarchical depth as a divisor invariant, and the explicit optimization theory developed over these canonical filtrations.

\medskip
\noindent
This places the present work in a complementary position relative to the existing literature on flag AG codes: the latter studies structural properties of nested code chains, whereas the present paper develops a geometric filtration mechanism that canonically generates and optimizes such chains.

The following example illustrates concretely the distinction between hierarchical AG codes and classical flag AG codes.

\begin{exam}
\label{ex:non-flag-distinction}
We illustrate here a feature of hierarchical AG codes that does not arise naturally
in the classical flag-code formalism.

Let $C/\F_{q^2}$ be the Hermitian curve
\[
C:\quad y^q+y=x^{q+1},
\]
with unique point at infinity $P_\infty$, genus
\[
g=\frac{q(q-1)}{2},
\]
and $n=q^3$ rational affine points.
Consider the canonical hierarchical filtration
\[
\L_0\subset \L_1\subset\cdots\subset \L_h,
\qquad
\L_i=\O_C(iP_\infty),
\]
induced by the hierarchical depth construction centered at $P_\infty$.

For each $i$, define the associated evaluation code
\[
\C_i=\C(C,iP_\infty,Z),
\qquad
Z=P_1+\cdots+P_n.
\]

This gives a nested chain
\[
\C_0\subset \C_1\subset\cdots\subset \C_h,
\]
which, as a chain of codes, is also a flag AG code in the classical sense.
However, the hierarchical structure contains additional canonical information:

\medskip
\noindent
\textbf{Step 1: Canonical determination of the filtration.}
In the flag-code formalism, one may choose many different increasing divisor sequences
\[
D_1\le D_2\le\cdots\le D_t.
\]
For instance, on the same Hermitian curve one may equally consider
\[
P_\infty,\; 2P_\infty,\; 4P_\infty,\; 7P_\infty,\dots
\]
or any other increasing chain.
Such choices are non-canonical.

By contrast, the hierarchical filtration above is uniquely determined by the divisor
$D=mP_\infty$ through the hierarchical depth invariant:
the sequence
\[
0,\; P_\infty,\; 2P_\infty,\;\dots,\; mP_\infty
\]
is intrinsic to the geometry of $D$ and not chosen ad hoc.

\medskip
\noindent
\textbf{Step 2: Canonical optimization index.}
Among all levels of the filtration, Proposition~\ref{prop:middle-layer} produces the
distinguished optimal index
\[
i^*=\left\lceil\frac{n+g-1}{2}\right\rfloor.
\]
Thus the filtration itself selects a preferred code
\[
\C_{i^*},
\]
determined canonically by the geometry of the curve and the divisor.

In the classical flag-code framework, one studies properties of a chosen nested family,
but there is generally no intrinsic mechanism selecting a unique optimal level from the flag.

\medskip
\noindent
\textbf{Conclusion.}
The nested chain of codes above may also be viewed as a flag AG code,
but its hierarchical interpretation provides two genuinely additional structures:

\begin{enumerate}
\item a canonical geometric origin of the filtration via hierarchical depth;
\item an intrinsic optimal-code selection principle through the explicit index $i^*$.
\end{enumerate}

These features go beyond the usual flag-code formalism, where nestedness alone does not
determine either a canonical filtration or a canonical optimal member.
\end{exam}

\subsection*{AG codes on surfaces}
\ \\ \\
Now we find out the optimal index $i$, which $\C_i$ is the optimal code in a nested sequence of codes obtained from hierarchical filtration of line bundles on a surface. Firstly, we reformulate Remark 2.2.5 from \cite{St} in term of AG codes corresponded to surfaces.

\begin{lem}
\label{lem:surface-des-dist-refined}
Let $S$ be a smooth projective surface over the finite field $\mathbb{F}_q$, and let $D$ be a very ample divisor on $S$. Define
\[
r := D^2 (q+1), \qquad d^* := |S(\mathbb{F}_q)| - r.
\]

Consider the evaluation code $\mathcal{C}(S,D,S(\mathbb{F}_q))$ obtained by evaluating global sections of $\mathcal{O}_S(D)$ at all rational points of $S$.

Then the minimum distance $d$ of $\mathcal{C}(S,D,S(\mathbb{F}_q))$ satisfies $d = d^*$ if and only if there exists a subset $Z' \subseteq S(\mathbb{F}_q)$ of cardinality $r$ such that
\[
h^0\big(S, \mathcal{I}_{Z'/S}(D)\big) > 0,
\]
where $\mathcal{I}_{Z'/S}$ denotes the ideal sheaf of $Z'$ in $S$.
\end{lem}

\begin{proof}
``Necessity.''
Assume that $d = d^* = |S(\mathbb{F}_q)| - r$. By definition of the minimum distance, there exists a nonzero section 
\[
s \in H^0\big(S, \mathcal{O}_S(D)\big)
\] 
whose evaluation vector
\[
(s(P))_{P \in S(\mathbb{F}_q)} \in \mathbb{F}_q^{|S(\mathbb{F}_q)|}
\] 
has weight $d^*$. This means that $s$ vanishes at exactly $r$ points of $S(\mathbb{F}_q)$. Let $Z' \subset S(\mathbb{F}_q)$ be the subset of these $r$ points. By construction, $s$ vanishes along $Z'$, so $s \in H^0\big(S, \mathcal{I}_{Z'/S}(D)\big)$, which implies
\[
h^0\big(S, \mathcal{I}_{Z'/S}(D)\big) > 0.
\]

``Sufficiency.''
Conversely, assume there exists a subset $Z' \subset S(\mathbb{F}_q)$ of cardinality $r$ such that 
\[
h^0\big(S, \mathcal{I}_{Z'/S}(D)\big) > 0.
\] 
Choose a nonzero section $s \in H^0\big(S, \mathcal{I}_{Z'/S}(D)\big)$. By construction, $s$ vanishes at all points of $Z'$, so the evaluation vector of $s$ has at least $r$ zero coordinates. Therefore, the weight of the evaluation vector satisfies
\[
\mathrm{wt}\big((s(P))_{P \in S(\mathbb{F}_q)}\big) \le |S(\mathbb{F}_q)| - r = d^*,
\]
so we conclude that $d \le d^*$.

Finally, applying Aubry's bound \cite[Proposition 3.1(ii)]{Au}, we have
\[
d \ge |S(\mathbb{F}_q)| - r = d^*,
\]
and therefore
\[
d = d^*.
\]
\end{proof}

Let $r$ be a positive real number. Define the function $\alpha:(0,r)\to\N\cap (0,r)$ with 
$$
\alpha(x)= \left\{
\begin{array}{rl}
	1&	0<x\leq 1	\\
	\lceil x\rfloor&	1\leq x\leq r-1 \\
	r-1&	r-1\leq x\leq r\ .
\end{array} \right.
$$

\begin{prop}
\label{prop:surface-optimal-i-refined}
Let $S$ be a smooth projective surface of arithmetic genus $g$ over $\mathbb{F}_q$. Denote $n := |S(\mathbb{F}_q)|$. Let $H$ be a very ample divisor on $S$ with $H^1(S,\O_S(iH))=H^2(S,\O_S(iH))=0$ for $i>0$. Set
$$a:=H^2>0,\quad b:=H.K_S,\quad c:=1+g,\quad r:=\sqrt{\frac{n}{a(q+1)}},$$
$D_i:=iH$ for all $i>0$ and $Z = S(\mathbb{F}_q)$.

For $0< i < r$, consider the evaluation codes
\[
\mathcal{C}_i := \mathcal{C}(S,D_i,Z),
\] 
whose dimension is
\[
k_i := \dim \mathcal{C}_i = \frac{1}{2} D_i \cdot (D_i - K_S) + 1 + g,
\]
and whose minimum distance is
\[
d_i := n - D_i^2 (q+1) = n - a(q+1)i^2.
\]

Denote by $i^*$ the index that maximizing 
\[
Q_i:= \frac{k_i}{n} d_i
\]
over $i \in (0,r)$. Then 

\begin{enumerate}
\item if $b\leq 0$ or $b\geq ar+\frac{2c}{r}$ we have
\[
i^* = \arg\max_{\substack{
i\in (0,r) \mathrm{\ is\ a\ real\ root\ of\ } \eqref{eq:surface-cubic-refined} \\
\mathrm{\ with}\ 6ai^2 - 3bi + 2c > ar^2
}} \{Q_{\alpha(i)}\}
\]

where
\begin{equation}
\label{eq:surface-cubic-refined}
4 a^2 (q+1) i^3 - 3 b a (q+1) i^2 + 2 a (-n + 2c(q+1)) i + b n = 0.
\end{equation}
But if for each real root of (\ref{eq:surface-cubic-refined}), $6ai^2-3bi+2c\leq ar^2$ then 
$$i^*= \arg \max\{Q_1,Q_{r-1}\}.$$

\item if $0< b< ar+\frac{2c}{r}$ then
\[
i^* = \arg\max_{\substack{
i\in (0,r) \mathrm{\ is\ a\ real\ root\ of\ } \eqref{eq:surface-cubic-refined}
}} \{Q_1,Q_{\alpha(i)}\}
\]
\end{enumerate}
 
\end{prop}

\begin{proof}
By the Riemann-Roch theorem for surfaces (\cite[p. 362]{Ha}), we have
\[
k_i = \dim H^0(S, \mathcal{O}_S(D_i)) = \frac{1}{2} D_i \cdot (D_i - K_S) + 1 + g 
= \frac{a}{2} i^2 - \frac{b}{2} i + c.
\]
%with $a = H^2$, $b = H \cdot K_S$, $c = 1+g$.  

By Lemma \ref{lem:surface-des-dist-refined} (with $Z = S(\mathbb{F}_q)$), the minimum distance satisfies
\[
d_i = n - D_i^2 (q+1) = n - a(q+1) i^2 > 0
\]
for $0<i <r$.

Define
\[
F(i) := k_i d_i = \left( \frac{a}{2} i^2 - \frac{b}{2} i + c \right) \left( n - a(q+1) i^2 \right).
\]

Maximizing $Q_i= F(i)/n$ over real $i$ is equivalent to maximizing $F(i)$. Differentiating:
\[
F'(i) = \left( a i - \frac{b}{2} \right) \left( n - a(q+1) i^2 \right) + \left( \frac{a}{2} i^2 - \frac{b}{2} i + c \right)( -2 a (q+1) i ).
\]

Simplifying and clearing factors of $\tfrac12$ yields the cubic equation
$$
4 a^2 (q+1) i^3 - 3 b a (q+1) i^2 + 2 a (-n + 2c(q+1)) i + b n = 0,
$$

which characterizes the critical points of $Q_i$ in the feasible interval $0<i<r$. 

We have 
\begin{center}
$F'(0)=-bn/2$\quad and\quad $F'(r)=2(bn-2acr(q+1)-anr)$.
\end{center}
But simplifying $F'(r)$ we obtain
\begin{align*}
F'(r) %& =-2(q+1)\big(2acr+\frac{n}{q+1}(ar-b)\big)\\
%&=-2(q+1)\big(2acr+ar^2(ar-b)\big)\\
%&=-2(q+1)ar\big(2c+r(ar-b)\big)\\
%&=-2(q+1)ar\big(2c+r(ar-b)\big)\\
&=-2(q+1)ar\big(ar^2-br+2c\big).
\end{align*}

%we obtain
%$$F'(r)=-2(q+1)ar(ar^2-br+2c).$$

It follows that $F'(0)F'(r)\leq 0$ if and only if $b(ar^2-br+2c)\leq 0$. By the assumptions (1) for $b$, $F'(0)F'(r)\leq 0$ and by the intermediate value theorem, the equation (\ref{eq:surface-cubic-refined}) has at least one real root in $(0,r)$. Suppose that $i$ is a real root of (\ref{eq:surface-cubic-refined}). If $$F''(i)=-6a^2(q+1)i^2+3ab(q+1)i-2ac(q+1)+an<0$$ then $i^*=i.$ But if for each real root $i$, $F''(i)\geq 0$ then $i^*= \arg \max\{Q_1,Q_{r-1}\}.$

When $0< b< ar+\frac{2c}{r}$ then $F'(0)<0$ and $F'(r)<0$. If $F'(t)\neq 0$ for all $t\in (0,r)$, then since $F'$ is continuous and it is negative at both ends of the interval, it follows from the intermediate value theorem (a continuous function that is negative at both endpoints and never zero stays negative everywhere in between) that $F$ is strictly decreasing on $(0,r)$. Therefore $i^*=1$. But if $F'(t)=0$ for some $t\in (0,r)$ then the maximum value for $F$ may occur in the critical points.
\end{proof}

\begin{exam}
Consider the projective plane $S=\P^2$ over $\F_7$. Let $H$ be the class of a line i.e., $H=\O_S(1)$ and let $Z$ is a reduced $\F_7$-rational zero-cycle of length $n=\# S(\F_7)=q^2+q+1=57$. Then $H^2=1$, $K_{\P^2}=\O_S(-3)$ and $g=0$. Thus $a=1$, $b=-3$ and $c=1$. Substituting these data in (\ref{eq:surface-cubic-refined}) gives the equation 
$$32i^3+72i^2-82i-171=0$$
with real roots $i\approx -2.38,\ -1.43,\ 1.56$. It follows from the range $0< i<\sqrt{57/8}\approx 2.66$ that $i^*=\lceil 1.56\rfloor=2$, i.e. the optimal code is $\C_2$. Of course, by calculating $Q_i$ for $i=1,2$, we obtain
$$Q_1\approx 2.57,\ Q_2\approx 2.63$$
which confirms the optimality of $i^*=2$.
\end{exam}

\begin{exam}
Let $S \subset \P^3_{\F_{q}}$ be a smooth quadric surface, with hyperplane class $H$. Let $H=F_1+F_2$ which $F_1$ and $F_2$ are rulings with $F^2_1=F^2_2=0$ and $F_1.F_2=1$. The number of $\F_q$-rational points on $S$ is $\# S(\F_q)=(q+1)^2.$ Let $q=13$. We have $H^2=\deg(S)=2$, $K_S=-2F_1-2F_2$, $H \cdot K_S =-4$, and
$$n=\# S(\F_{13})=(13+1)^2=196.$$
For $D=iH$, Riemann-Roch (with $p_a(S)=0$) yields
$$\dim \C(S,iH,Z)=\ell(iH)=\frac{1}{2}\left(i^2 H^2-iH\cdot K_S \right)+1=i^2+2i+1.$$

The performance functional is
$$Q_i=\frac{\ell(iH)}{n}\cdot d_i=\frac{i^2+2i+1}{196}\cdot\left(196-28i^2\right).$$
The admissible range $d_i>0$ implies $i\leq 2$. The equation (\ref{eq:surface-cubic-refined}) gives the real roots $i\approx -2.13,\ -1,\ 1.63$. So $i^*=\lceil 1.63\rfloor=2$.
Also, a direct check shows
$$Q_1\approx 3.43,\quad Q_2\approx 3.86.$$
%Thus the optimal index is $i^*=2$.

\end{exam}

\begin{rem}
Note that in last two examples, the value of $\sqrt{n/(H^2(q+1))}$ is small, so the candidate set $\{i\}$ is very short ($1,2$ only). This makes the optimum highly predictable and preserves high minimum distance, but limits the diversity of available code dimensions. In the next example, we will solve this problem to some extent.
\end{rem}

\begin{exam}
Let $S\subset\P^3_{\F_{q^2}}$ be the Hermitian surface
$$X_0^{q+1}+X_1^{q+1}+X_2^{q+1}+X_3^{q+1}=0.$$
It is smooth of degree $d = q+1$.
Let $H$ be the hyperplane class. Then
$$H^2=q+1,\quad K_S\sim (q-3)H,\quad H\cdot K_S=(q-3)(q+1),$$
and
$$n=\#S(\F_{q^2})=(q^3+1)(q^2+1).$$
See, for example, \cite{HoKi} for details about Hermitian surfaces.

(1) Let $D=iH$ with $i>q-3$. Then we have $H^1(S,\O_S(iH))=H^2(S,\O_S(iH))=0$. Hence the assumptions of Proposition \ref{prop:surface-optimal-i-refined} holds. The Riemann-Roch (with arithmetic genus $g(S)=\binom{q}{3}$) yields
$$\ell(iH)=\frac{q+1}{2}\left(i^2-i(q-3)\right)+1.$$

Thus the performance functional is
$$Q_i=\frac{\ell(iH)}{n}\cdot d_i=\frac{\frac{q+1}{2}\left(i^2-i(q-3)\right)+1}{n}\cdot \left(n-i^2(q+1)^2\right).$$

For $q=5$, we have $H^2=6$, $H\cdot K_S =12$, $n=3,276$,
$$\ell(iH)=3(i^2-2i)+1,\quad d_i\geq 3,276-156 i^2.$$
The admissible range $d_i>0$, $i>q-3=2$ and $i<\sqrt{n/(H^2(q+1))}\approx 9.53$ give the range $3\leq i\leq 9$. Moreover, the equation (\ref{eq:surface-cubic-refined}) has the real roots 
$$-7.82,\ 1.11,\ 5.20.$$ 
Thus the optimum index is $i^*=5$ and, moreover, $Q_5\approx 33.36$.

For $q=101$, we have $H^2=102$, $H\cdot K_S =9,996$, $n=10,511,141,004$ and, also, by $i>q-3=98$ and $i<\sqrt{n/(H^2(q+1))}\approx 1,005.13$ we obtain the range $99\leq i\leq 1,005$. Substituting in the equation (\ref{eq:surface-cubic-refined}) gives the real roots $$-698.84,\ 48.96,\ 723.38.$$ It follows that the optimum index is $i^*=723$ and $Q_{723}\approx 11,119,916$.

(2) Let $D=iH$ with $i\leq q-3$. Then we have $H^1(S,\O_S(iH))=0$ and

$$\h^2(S,O_S(iH))=\h^0(S,O_S((q-3-i)H)=\h^0(S,O_{P^3}((q-3-i))=\binom{q-i}{3}.$$

On the other hand, 
$$k_i=\h^0(S,O_S(iH))=\frac{1}{2}D\cdot (D-K_S)+\chi(\O_S)+\h^1(S,O_S(D))-\h^2(S,O_S(iH))$$

and $\chi(\O_S)=1+\frac{q(q-1)(q-2)}{6}$, so
$$k_i=\frac{1}{2}i(i-q+3)(q+1)+1+\frac{q(q-1)(q-2)}{6}-\binom{q-i}{3}.$$
Also, we have $d_i=n-i^2(q+1)^2$. These data follow that
%$$Q_i=\frac{1}{n}\big(\frac{1}{2}i(i-q+3)(q+1)+1+\frac{q(q-1)(q-2)}{6}-\binom{q-i}{3}).(n-i^2(q+1)^2)$$

\[
Q_i =\frac{1}{6n}\big( 6n 
+ 11n i
+ 6(n - (q + 1)^2) i^2
+ (n - 11 (q + 1)^2) i^3
- 6(q + 1)^2 i^4
- (q + 1)^2 i^5\big).
\]
Set $f(i):=6nQ_i$. Then 
$$f'(i)=11n+12(n - (q + 1)^2)i+3(n - 11 (q + 1)^2)i^2-24(q + 1)^2 i^3-5(q + 1)^2 i^4.$$

For $q=5$ and $i=1,2$ we have $Q_1\approx 3.96$ and $Q_2\approx 9.56$. %Therefore the optimal index $i^*=5$ obtained above is approved.

For $q=101$ we have $n=10,511,141,004$ and the roots of $f'(i)=0$ are, approximately,
$$-778.97,\ -2.57,\ -1.42,\ 778.17.$$
Since $i\leq q-3=98$, so none of these roots are confirmed. 

Consequently, the optimal indices in the cases $q=5$ and $q=101$ are $5$ and $723$, respectively.
\end{exam}

\begin{rem}
For Hermitian surfaces, the bound
$$i < \sqrt{\frac{n}{H^2(q+1)}} = \sqrt{\frac{(q^3+1)(q^2+1)}{(q+1)^2}}$$
controls the range of $i$. For moderate $q$, this range is significantly larger than in the quadric surface case, allowing a richer choice of code parameters, though at the cost of lower minimum distance for large $i$.
\end{rem}

\section*{Data Availability Statement.} 
No datasets were generated or analyzed during the current study.

\def\baselinestretch{1}

\def\baselinestretch{1.66}

\ \\ \\
Rahim Rahmati-Asghar,\\
Department of Mathematics, Faculty of Basic Sciences,\\
University of Maragheh, P. O. Box 55181-83111, Maragheh, Iran.\\
E-mail:  \email{rahmatiasghar.r@maragheh.ac.ir, rahmatiasghar.r@gmail.com}

\end{document}